\documentclass[openany]{article}
\usepackage[margin=1in]{geometry}
\usepackage{amsmath}

\newcommand{\abs}[1]{\left|#1\right|}

\newenvironment{nproof}[1]{\noindent\textit{Proof{#1}.}}{\hfill Q.E.D.}

\newtheorem{theorem}{Theorem}[section]
\newtheorem{cor}[theorem]{Corollary}

\newtheorem{rem}[theorem]{Remark}

\newtheorem{lem}[theorem]{Lemma}
\newtheorem{example}[theorem]{Example}

\DeclareMathOperator{\Hess}{Hess}

\DeclareMathOperator{\sign}{sign}

\title{On the Pontrjagin classes of spray manifolds}
\author{ Zhongmin Shen\ and \   Runzhong Zhao}
\date{}
\begin{document}
\maketitle

\begin{abstract}

The characterization of projectively flat Finsler  metrics on an open subset in $R^n$ is the Hilbert’s Fourth Problem in the regular case. Locally projective flat Finsler manifolds form an important class of Finsler manifolds. Every Finsler metric induces a spray on the manifold via geodesics. Therefore it is a natural problem to investigate the geometric and topological properties of  manifolds equipped with a spray. 
In this paper we study the Pontrjagin classes of a manifold equipped with a locally projectively flat spray and show that such manifold must have zero Pontrjagin classes.
\end{abstract}

\section{Introduction}

The notion of projectively flat metrics in Finsler geometry originated from the famous Hilbert's Fourth Problem, that is, to characterize Finsler metrics on an open subset of $\mathbf{R}^n$ whose
geodesics are straight lines as sets of points. 
 By  Beltrami's theorem, a Riemannian metric is locally
projectively flat if and only if it is of constant sectional curvature. However, things are much more complicated in Finsler  geometry. Although every locally projectively flat
Finsler metric is of scalar flag curvature, many are not of constant flag curvature\cite{HuMo}\cite{Zhou}. There are also Finsler  metrics of scalar flag curvature but are not locally projectively flat.
Therefore  the study of locally projectively flat metrics is of its own interest. The characterization of locally projectively flat Finsler metrics dates
back to Hamel's work\cite{Ha} at the beginning of 20th century. In recent years, many more explicit constructions have been found for special types of Finsler metrics.
For example, a Randers metric $F = \alpha + \beta$ is locally projectively flat if and only if $\alpha$ is a locally projectively flat Riemannian metric and $\beta$ is a closed
one-form\cite{BaMa}. Locally projectively flat $(\alpha, \beta)$ metrics and general $(\alpha,\beta)$ metrics are also studied intensively (see, for example, \cite{LiSh}\cite{Sh2}\cite{Zhu}).

Every Finsler metric induces a spray via geodesics.  The study of sprays on a manifold will lead to a better understanding on Finsler metrics. The notion of projective flatness can be extended naturally to sprays.  We would like to know any topological obstruction to the existence of locally projectively flat sprays on a manifold.

\begin{theorem}\label{MTM}
	Let $G$ be a locally projectively flat spray on a $n$-manifold $M$. Then the Pontrjagin classes $p_k$ of $M$ with rational coefficients are all zero for $k\ne 0$.
\end{theorem}
This gives a topological obstruction to the existence of locally projectively flat sprays on a manifold.
\begin{example}
	It is well known that spheres have trivial Pontryagin classes. The standard Riemannian metric $g_{\mathbf{S}^n}$ on an $n$-dimensional unit sphere has constant sectional curvature $1$, 
	hence is locally projectively 
	flat. Let $G_0$ be the spray of this metric, a family of locally projectively flat sprays can be constructed using positively homogeneous functions $P$ on the slit tangent bundle 
	$T\mathbf{S}^n\setminus 0$. See section \ref{EX} for more detailed constructions.
\end{example}
\begin{example}
	The cohomology of the complex projective space $\mathbf{C}P^n$ with coefficient $G$ is
	\[
		H^l(\mathbf{C}P^n,G) = \begin{cases} G & \text{if\ }l\text{\ is even and }0 \le l \le n \\ 0 & \text{otherwise}\end{cases}
	\]
	and the total Pontrjagin class of $\mathbf{C}P^n$ is $p = (1+c^2)^{n+1}$ where $c$ is the image of the generator of $H^2(\mathbf{C}P^n,\mathbf{Z})$ under the natural map
	$H^2(\mathbf{C}P^n,\mathbf{Z}) \to H^2(\mathbf{C}P^n,G)$. In particular, the Pontrjagin classes of $\mathbf{C}P^n(n\ge 2)$ are not zero, hence $\mathbf{C}P^n(n\ge 2)$ does not admit
	a locally projectively flat spray.
\end{example}

\section{Preliminaries}

A spray $G$ on a manifold $M$ is a smooth vector field on the slit tangent bundle $TM\setminus 0$ expressed in a standard
local coordinate system $(x^i, y^i)$ on $TM$ as
\begin{equation}
	G = y^i \frac\partial{\partial x^i} - 2G^i(x,y)\frac\partial{\partial y^i}
\end{equation}
where $G^i(y)$ are local functions on $TM$ satisfying
\[
	G^i(x,\lambda y) = \lambda^2 G^i(x,y)
\]
for all $\lambda > 0$. The $G^i$'s are also called the \textit{spray coefficients}.

Let $\tilde{\gamma}$ be an integral curve of $G$ on $TM\setminus 0$, and $\gamma = \pi\circ\tilde{\gamma}$ be its projection
on $M$, where $\pi: TM\setminus 0\to M$ is the canonical projection. $\gamma$ satisfies the equation
\[
	\ddot{\gamma}^i + 2G^i(\gamma, \dot{\gamma}) = 0
\]
and it is called a \textit{geodesic} of $G$. We say that two sprays $G$ and $\tilde{G}$ are \textit{(pointwise) projectively related} if their geodesics are the same as sets of points on the manifold.
Equivalently, this is characterized by the condition
\[
	\tilde{G}^i = G^i + Py^i
\]
where $P = P(x,y)$ satisfies the homogeneity property
\[
	P(x,\lambda y) = \lambda P(x,y)
\]
for $\lambda > 0$. A spray $G$ is said to be \textit{flat} if at every point, there is a local coordinate system in which
\[
	G = y^i\frac{\partial}{\partial x^i}.
\]
It is said to be \textit{locally projectively flat} if it is projectively related to a flat spray.
A quantity is said to be \textit{projectively invariant} if it is the same for projectively related sprays.

In order to study the Pontrjagin classes of $\pi^*TM$ and $TM$ we shall use the \textit{Berwald connection} on the
pullback bundle $\pi^*TM$. Let
\[
	N^i_j := \frac{\partial G^i}{\partial y^j}, \quad\Gamma^i_{jk} := \frac{\partial^2 G^i}{\partial y^j\partial y^k},
\]
the local connection 1-forms of the Berwald connection are given by
\begin{equation}
	\omega_j^{\;i} = \Gamma^i_{jk}dx^k
\end{equation}
and its curvature forms are
\begin{equation}\label{curvform}
	\Omega_j^{\;i} := d\omega_j^{\;i} - \omega_j^{\;k}\wedge \omega_k^{\;i}.
\end{equation}
Put $\omega^i := dx^i$ and $\omega^{n+i} := dy^i + N^i_jdx^j$ we have
\[
	d\omega^i = \omega^j \wedge \omega_j^{\;i}
\]
and
\begin{equation}
	\Omega_j^{\;i} = \frac12R_{j\;kl}^{\;i}\omega^k\wedge \omega^l - B_{j\;kl}^{\;i}\omega^k\wedge\omega^{n+l}
\end{equation}
where $R$ is the \textit{Riemann curvature} and $B$ is the \textit{Berwald curvature}. In local coordinates, they are given by
\[\begin{aligned}
	&R_{j\;kl}^{\;i} = \frac{\delta\Gamma^i_{jl}}{\delta x^k} - \frac{\delta\Gamma^i_{jk}}{\delta x^l}
	+ \Gamma^i_{km}\Gamma^m_{jl} - \Gamma^i_{lm}\Gamma^m_{jk} \\
	&B_{j\;kl}^{\;i} = \frac{\partial^3G^i}{\partial y^j\partial y^k\partial y^l}
\end{aligned}\]
where $\frac{\delta}{\delta x^i} = \frac{\partial}{\partial x^i} - N^j_i\frac{\partial}{\partial y^j}$ is the horizontal covariant derivative. For the simplicity of notation we will denote
$A^*_{\;*|j} = \frac{\delta}{\delta x^j}A^*_{\;*}$ and $A^*_{\;*\cdot j} = \frac{\partial}{\partial y^j}A^*_{\;*}$.

The two-index Riemannian curvature tensor is then given by $R^i_{\;k} = R^{\;i}_{j\;kl}y^jy^l$. We have
\begin{equation}
	R^{\;i}_{j\;kl} = \frac13\left(R^i_{\;k\cdot l\cdot j} - R^i_{\;l\cdot k\cdot j}\right)
\end{equation}
so the two-index Riemann curvature tensor and the four-index Riemann curvature tensor
basically contain the same geometric data.
The Riemann curvature can be computed directly using the spray coefficients as
\begin{equation}\label{RC}
	R^i_{\;k} = 2\frac{\partial G^i}{\partial x^k} - y^j\frac{\partial^2G^i}{\partial x^j\partial y^k} + 2G^j\frac{\partial^2G^i}{\partial y^j\partial y^k} - \frac{\partial
		G^i}{\partial y^j}\frac{\partial G^j}{\partial y^k}.
\end{equation}
We will also use $R = \frac1{n-1}R^m_{\;m}$.

In the case when  $G$ is a Berwald spray, where $B_{j\;kl}^{\;i} = 0$, the spray coefficients $G^i$'s are quadratic in $y^i$'s, and so are the Riemann curvature $R^i_{\;k}$. In particular,
$\Gamma^i_{jk}$ and $R^{\;i}_{j\;kl}$ are independent of $y^i$'s. Therefore, the differential forms $\omega^{\;i}_j$ and $\Omega^{\;i}_j$ can be viewed as differential forms on
$M$, and the connection can be viewed as a connection of the tangent bundle $TM$. In fact, we have
\[
	\Omega^{\;i}_j = \dfrac12R_{j\;kl}^{\;i}dx^k\wedge dx^l.
\]

In dimension  $n \geq 3$, locally projectively flat sprays are characterized by two projectively invariant quantities.
The \textit{Douglas curvature} is constructed from the Berwald curvature:
\begin{equation}\label{Douglas}
	D_{j\;kl}^{\;i} = B_{j\;kl}^{\;i} - \frac2{n+1}\left[
	E_{jk}\delta^i_l + E_{jl}\delta^i_k + E_{kl}\delta^i_j + \frac{\partial E_{jk}}{\partial y^l}y^i \right]
\end{equation}
where $E_{ij} =\frac{1}{2} B^{\;m}_{m\;ij}$ is called the \textit{mean Berwald curvature}.

The \textit{Weyl tensor} is the spray analog of the projective curvature tensor in Riemannian geometry. It is defined by
\begin{equation}
	W^i_{\;k} = A^i_{\;k} - \frac1{n+1}\frac{\partial A^m_{\;k}}{\partial y^m}y^i
\end{equation}
where $A^i_{\;k} = R^i_{\;k} - R\delta^i_k$.
The followings are well-known\cite{Sh3}:
\begin{lem}
	A spray $G$ is of scalar curvature, in the sense that
	\[
		R^i_{\;k} = R\delta^i_k - \tau_ky^i
	\]
	for some 1-form $\tau$ on $TM\setminus 0$ with $\tau_k y^k = R$, if and only if $W = 0$.
\end{lem}

\begin{lem}
	A spray $G$ on a manifold $M$ of dimension $\ge 3$ is locally projectively flat if and only if $W = 0$ and $D = 0$.
\end{lem}

For convenience in discussing the projective change by the $S$-curvature we will also use another formula for the Weyl tensor.
Let $dV = \sigma(x)dx^1\wedge\cdots\wedge dx^n$ be a volume form on $M$, the quantity
\[
	\tau(x,y) := \ln \frac{\sqrt{\det g_{ij}(x,y)}}{\sigma(x)}
\]
is called the \textit{distortion} and its rate of change along geodesics is measured by \textit{$S$-curvature}.
Namely, let $\gamma(t)$ be a geodesic with $\gamma(0) = x$ and $\dot{\gamma}(0) = y\in T_xM\setminus 0$, we have
\begin{equation}
		\mathbf{S}(x,y) := \left.\frac{d}{dt}\right|_{t=0}\left[\tau(\gamma(t), \dot{\gamma}(t))\right]
\end{equation}
The $S$-curvature can be expressed as
\[
	S(x,y) = \frac{\partial G^m}{\partial y^m}(x,y) - y^m\frac{\partial}{\partial x^m}[\ln \sigma ](x)
\]
and satisfies the homogeneity property
\[
	S(x,\lambda y) = \lambda S(x, y)
\]
for $\lambda > 0$.
The non-Riemannian quantity \textit{$\chi$-curvature}
\[
	\chi_k := -\frac16\left(2R^m_{\;k\cdot m} + R^m_{\;m\cdot k}\right).
\]
relates the $S$-curvature and the Weyl tensor. Indeed, it can be computed using the $S$-curvature by
\[
	\chi_k = \frac12\left(S_{\cdot k|m}y^m - S_{|k}\right)
\]
and the Weyl tensor can be expressed as
\[
	W^i_{\;k} = R^i_{\;k} - \left(R\delta^i_k - \frac12R_{\cdot k}y^i\right) + \frac3{n+1}\chi_ky^i.
\]
Thus for a spray $(G, dV)$ whose $S$-curvature is vanishing, we have $\chi_k = 0$ and
\begin{equation}\label{Weyl}
	W^i_{\;k} = R^i_{\;k} - \left(R\delta^i_k - \frac12R_{\cdot k}y^i\right).
\end{equation}

\section{Pontrjagin Classes}\label{section3}

The Pontrjagin classes were originally introduced by Lev Pontrjagin in 1940's in the study of the Grassmannian manifolds. In modern texts, it is most often described using the Chern class of a
complexified bundle.
Let $\xi$ be a real vector bundle, and $\xi_\mathbf{C} := \xi\otimes_\mathbf{R}\mathbf{C}$ be its complexification, then the $i$-th Pontrjagin class of $\xi$ is given by
\[
	p_i = (-1)^ic_{2i}(\xi_\mathbf{C})
\]
where $c_{j}(\xi_\mathbf{C})$ is the $j$-th Chern class of the complex vector bundle $\xi_\mathbf{C}$.

Using a connection $D$ on the vector bundle $\xi$, we can express the Pontrjagin classes of $\xi$ as a differential form representing a de Rham cohomology class.
Let $U\subset M$ be an open subset on which there is trivialization $\psi: \xi|_U \cong U\times \mathbf{R}^n$ of the bundle $\xi$, and $e_i = \psi^{-1}\partial_i$.
The connection $D$, viewed as a covariant derivative, can be described by
\[
	D_X  \eta  := \left[d\eta^i(X) + \eta^j\omega^{\;i}_j(X)\right]e_i
\]
where $X$ is a vector field on $U$, $\eta$ is a section of $\xi$ given by $\eta = \eta^ie_i$ on $U$, and the $\omega_j^{\;i}$'s are $n^2$ local one-forms on $U$.
The curvature of $D$ is hence described by an $n\times n$ matrix of (real-valued) 2-forms
\[
	\Omega = (\Omega_j^{\;i})
\]
where $\Omega_j^{\;i}$'s are given as in (\ref{curvform}).
The differential forms $\omega_j^{\;i}$'s and $\Omega_j^{\;i}$'s, viewed as complex-valued forms, define a connection on the complexified bundle $\xi_\mathbf{C}$ and its curvature, which we will still
denote by $D$ and $\Omega$ respectively.

Given an $n\times n$ complex matrix $T$, we have
\begin{equation}\label{det}
	\det(I+tT) = 1 + t\sigma_1(T)+ \cdots + t^n\sigma_n(T)
\end{equation}
where $\sigma_i(T)$ is the $i$-th elementary symmetric polynomial evaluated on the $n$ eigenvalues of $T$. The following can be found on most textbooks covering characteristic classes (see, for
example, \cite{MiSt}):

\begin{theorem}
	Let $\xi$ be a complex vector bundle of rank $n$ with connection $D$. Then the cohomology class $[\sigma_r(\Omega)]\in H^{2r}(X;\mathbf{C})$
	is equal to $(2\pi i)^rc_r(\xi)$, for all $r = 1, \ldots, n$, where $c_r(\xi)$ is the $r$-th Chern class with coefficient $\mathbf{C}$.
\end{theorem}

In the study of Riemannian manifolds of constant curvature, Chern obtained a consequence of this theorem\cite{Ch1}\cite{MiSt}:

\begin{cor}
	Let $\xi$ be a real vector bundle of rank $n$ with connection $D$. Then the de Rham cocycle $\sigma_{2k}(\Omega)$ represents the cohomology class $(2\pi)^{2k}p_k(\xi)\in H^{4k}(M;\mathbf{R})$, while
	$[\sigma_{2k+1}(\Omega)] = 0$ in $H^{4k+2}(M;\mathbf{R})$, where $p_k(\xi)$ is the $k$-th Pontrjagin class of $\xi$ with coefficient $\mathbf{R}$.
\end{cor}

Chern proved the following

\begin{theorem}
	Suppose that for a Riemannian manifold $M$, the sectional curvature $K(V_x, W_x)$ of the plane $\mathop{\rm span}\{V_x, W_x\}\subset T_xM$ depends on the point $x\in M$ only. Then all its Pontrjagin
	classes with rational coeffecients are zero.
\end{theorem}

We shall remark that such manifolds are said to have \textit{isotropic curvature} in Finsler geometry.

In terms of the local curvature forms $\Omega_j^{\;i}$'s defined in some neighborhood of $x\in M$, a direct computation using (\ref{det}) gives
\[\begin{aligned}
	\sigma_{2k}(\Omega)(x) =& \sum_{i_1 < \cdots < i_{2k}, \sigma\in S_{2k}} \sign{\sigma}\cdot \Omega_{i_1}^{\;i_{\sigma(1)}}(x) \wedge \cdots\wedge
	\Omega_{i_{2k}}^{\;i_{\sigma(2k)}}(x) \\
	=& \dfrac1{(2k)!}\sum_{i_1, \ldots, i_{2k}, \sigma\in S_{2k}}\sign{\sigma}\cdot \Omega_{i_1}^{\;i_{\sigma(1)}}(x) \wedge \cdots\wedge
	\Omega_{i_{2k}}^{\;i_{\sigma(2k)}}(x)
\end{aligned}\]
where $S_{2k}$ is the symmetric group on $2k$ elements.


\section{Pontrjagin Classes of Douglas Sprays}

In this section we study the Pontrjagin classes of Douglas sprays. The main tool will be a projective change by the $S$-curvature.

Let $G$ be a spray and $dV$ be a volume form on an $n$-dimensional manifold $M$.
We define another spray $\hat{G}$ by
\begin{equation}
	\hat{G}^i = G^i - \frac{S}{n+1}y^i
\end{equation}
where $S$ is the $S$-curvature of $(G, dV)$. In the sequel a letter with a hat over it will always represent a quantity of the spray $\hat{G}$.
We have the following
\begin{lem}{\cite{Sh4}}
	\begin{equation}\begin{aligned}
		\hat{B}_{j\;kl}^{\;i} =& D_{j\;kl}^{\;i}\\
		\hat{S} =& 0
	\end{aligned}\end{equation}
\end{lem}

Thus if $G$ is Douglas, then $\hat{G}$ is a Berwald spray.
As an immediate consequence, we have
\begin{lem}
	Let $G$ be a Douglas spray on a manifold $M$. The Pontrjagin classes of $M$ with coefficient $\mathbf{R}$ are represented by the forms
	\[
		\frac1{(4\pi)^{2k}(2k)!}\sum_{i_1, \ldots ,i_{2k}, \sigma\in S_{2k}}
		\sign \sigma\cdot \prod_{s=1}^{2k}\hat{R}^{\;i_{\sigma(s)}}_{i_s\;m_sl_l}\bigwedge_{s=1}^{2k}dx^{m_s}\wedge dx^{l_s}
	\]
	where $\hat{R}$ is the Riemann curvature of the associated spray $\hat{G}$.
\end{lem}

\begin{nproof}{ of Theorem \ref{MTM}}
	We will assume that the dimension of the manifold is $n\ge 4$, for otherwise there is nothing to prove.
	We already have that $\hat{G}$ is a Berwald spray. On the other hand, $\hat{S} = 0$ implies that $\hat{\chi}_k = 0$ for the spray $\hat{G}$.
	The expression for the Weyl tensor (\ref{Weyl}) now becomes
	\[
		\hat{W}^i_{\;k} = \hat{R}^i_{\;k} - \left(\hat{R}\delta^i_k - \frac12\hat{R}_{\cdot k}y^i\right) = 0
	\]
	hence
	\[
		\hat{R}^i_{\;k} = \left(\hat{R}\delta^i_k - \frac12\hat{R}_{\cdot k}y^i\right).
	\]
	A straightforward calculation yields
	\[
		\hat{R}^{\;i}_{j\;kl} = \frac13\left(\hat{R}^i_{\;k\cdot l\cdot j} - \hat{R}^i_{\;l\cdot k\cdot j}\right) = \frac12\left(\hat{R}_{\cdot l\cdot j}\delta^i_k - \hat{R}_{\cdot k\cdot j}\delta^i_l
		\right).
	\]
	It follows that the curvature forms of the Berwald connection of $\hat{G}$ is given by
	\[
		\hat{\Omega}_j^{\;i} = \frac14\left(\hat{R}_{\cdot l\cdot j}\delta^i_k - \hat{R}_{\cdot k\cdot j}\delta^i_l\right)dx^k\wedge dx^l  = \frac{1}{2} \hat{R}_{\cdot l\cdot j}\; dx^i\wedge dx^l 
	\]
Then
\begin{eqnarray*}
		\sigma_{2k}(\hat{\Omega}) &=& \frac1{(2k!)}\sum_{i_1, \ldots ,i_{2k}, \sigma\in S_{2k}} \sign{\sigma}\cdot \hat{\Omega}_{i_1}^{\;i_{\sigma(1)}} \wedge \cdots\wedge
		\hat{\Omega}_{i_{2k}}^{\;i_{\sigma(2k)}} \\
& = & \frac{1}{2^{2k}(2k)!} \sum_{i_1, \ldots ,i_{2k}, \sigma\in S_{2k}} \sign{\sigma} \cdot \hat{R}_{\cdot l_1\cdot i_1} \cdots \hat{R}_{\cdot l_{2k}\cdot i_{2k}} 
dx^{i_{\sigma(1)}}\wedge dx^{l_1} \cdots dx^{i_{\sigma(2k)}}\wedge dx^{l_{2k}}\\
& = & \frac{(-1)^{k}}{2^{2k}(2k)!} \sum_{i_1, \ldots ,i_{2k}, \sigma\in S_{2k}} \sign\sigma  \cdot \hat{R}_{\cdot l_1\cdot i_1} \cdots \hat{R}_{\cdot l_{2k}\cdot i_{2k}} 
dx^{i_{\sigma(1)}}\wedge\cdots \wedge  dx^{i_{\sigma(2k)}}\wedge dx^{l_1}\wedge \cdots \wedge  dx^{l_{2k}}\\
& = & \frac{(-1)^{k}}{2^{2k}} \sum_{i_1, \ldots ,i_{2k}}  \hat{R}_{\cdot l_1\cdot i_1} \cdots \hat{R}_{\cdot l_{2k}\cdot i_{2k}} dx^{i_{1}}\wedge\cdots \wedge  dx^{i_{2k}} \wedge
dx^{l_1}\wedge \cdots \wedge  dx^{l_{2k}}
\end{eqnarray*}
	Since $\hat{R}_{\cdot i\cdot j} = \hat{R}_{\cdot j\cdot i}$, we obtain $0$ after summing over $i_1$ and $l_1$. We conclude that
	\[
		\sigma_{2k}(\hat{\Omega}) = 0.
	\]
	Since $\hat{G}$ is Berwald, these curvature forms can be viewed as curvature forms on $TM$. Therefore the Pontrjagin classes of $M$ with real coefficients are identically zero.
	It follows from naturality of cohomology in coefficients that the Pontrjagin classes of $M$ with rational coefficients are zero.
\end{nproof}

\section{Locally Projectively Flat Sprays on the Sphere}\label{EX}

	In this section we describe in more detail some known locally projectively flat sprays on the unit sphere.
	Let $G_0$ be the spray of the standard Riemannian metric on the unit sphere, and $P: T\mathbf{S}^n\setminus 0\to\mathbf{R}$ be a smooth function satisfying 
	\[
		P(x, \lambda y) = \lambda P(x,y)
	\]
	for all $x\in \mathbf{S}^n$, $y\in T_x\mathbf{S}^n\setminus\{0\}$, and $\lambda > 0$. Then the spray $G_P = G_0 - 2PY$ is locally projectively flat, where $Y$ is the vector field expressed as 
	$y^i\frac{\partial}{\partial x^i}$ in local coordinates $(x^i, y^i)$ for the slit tangent bundle $T\mathbf{S}^n\setminus 0$. 

	Among this family are a collection of sprays induced by Randers metrics $F(y) = \sqrt{g_{\mathbf{S}^n}(y,y)} + df(y)$ 
	where $f$ is a smooth function on $\mathbf{S}^n$ and $y\in T\mathbf{S}^n\setminus 0$. In this case we have 
	\[
		P(y) = \dfrac{\Hess f(y,y)}{2F(y)}.
	\]

	Another interesting class of locally projectively flat Finsler metrics on the sphere are the Bryant metrics\cite{Br}\cite{Sh1}. In gnomonic coordinates it is given by 
	\[
		F(y) = \sqrt{\dfrac{\sqrt{A}+B}{2D}+\left(\dfrac{C}{D}\right)^2} + \dfrac{C}{D}
	\]
	where 
	\[\begin{aligned}
		A =& \left(\cos (2\alpha) \abs{y}^2 +\left(\abs{x}^2\abs{y}^2 - \langle x,y\rangle^2\right)\right)^2 + \left(\sin (2\alpha) \abs{y}^2\right)^2 \\
		B =& \cos (2\alpha) \abs{y}^2 +\left(\abs{x}^2\abs{y}^2 - \langle x,y\rangle^2\right) \\
		C =& \sin (2\alpha) \langle x,y\rangle \\
		D =& \abs{x}^4 + 2\cos (2\alpha) \abs{x}^2 + 1
	\end{aligned}\]
	with $\langle\cdot,\cdot\rangle$ being the standard inner product in the Euclidean space $\mathbf{R}^n$ and $\abs{\cdot}$ being the induced norm.
	The spray induced by a Bryant metric can also be written in the form $G_0 - 2PY$, with the function $P$ satisfying
	\[
		\left[(P + dp(y))^2 + (dq(y)^2)\right]\left[(F(y) - dq(y))^2 - dq(y)^2\right] = (\Hess r(y,y) + s(x)g_{\mathbf{S}^n}(y,y))^2
	\]
	where $y\in T\mathbf{S}^n\setminus 0$ and the functions $p,q,s$ can be expressed in the gnomonic coordinates as
	\[\begin{aligned}
		p(x) =& \dfrac14\ln \dfrac{1+2\cos(2\alpha) \abs{x}^2 + \abs{x}^4}{1 + 2\abs{x}^2 + \abs{x}^4} \\
		q(x) =& \dfrac1{2}\tan^{-1}\dfrac{\abs{x}^2 + \cos(2\alpha)}{\sin(2\alpha)} \\
		s(x) =& s(\abs{x}^2) = 4(1+\abs{x}^2)^2r''(\abs{x}^2) + 4(1+\abs{x}^2)r'(\abs{x}^2)
	\end{aligned}\]
	and $r = r(\abs{x}^2)$ satisfies the differential equation 
	\begin{equation}\label{DEP}
		4r''(\abs{x}^2)(1+\abs{x}^2) + 6r'(\abs{x}^2) = \dfrac{\sin(2\alpha)}{2(1+2\cos(2\alpha)\abs{x}^2 + \abs{x}^4)}.
	\end{equation}
	One may check that these functions are indeed smooth on the whole sphere as follows.
	Let $\varphi$ be the inclination angle in the spherical coordinate. Since the Bryant metrics are invariant under vertical reflection $\varphi\mapsto \pi - \varphi$, it suffices to check that these 
	functions extend to even smooth functions in $u = 1/\tan\varphi = 1/\abs{x}$. For $p$ and $q$ this is clear. For $s$ we have 
	\[
		s(u) = u^2(1+u^2)^2\dfrac{d^2r}{du^2} + (3u^3 + 5u)(u^2 + 1)\dfrac{dr}{du}
	\]
	so $s$ is even and smooth as long as $r$ is. Finally, the differential equation (\ref{DEP}) reads 
	\[
		(1+u^2)\dfrac{d^2r}{du^2} + 3u\dfrac{dr}{du} = \dfrac{\sin (2\alpha)}{2(1 + 2\cos(2\alpha)u^2 + u^4)}.
	\]
	Note that it can be written as 
	\[
		\dfrac{d}{du}\begin{pmatrix} \dot{r} \\ r \end{pmatrix} = f(u,r,\dot{r}) := \begin{pmatrix} \dfrac{3u}{1+u^2}\dot{r} \\ \dot{r} \end{pmatrix} + 
		\begin{pmatrix} \dfrac{\sin (2\alpha)}{2(1 + 2\cos(2\alpha)u^2 + u^4)(1+u^2)} \\ 0 \end{pmatrix}
	\]
	Since $f$ is uniformly Lipschitz in $(r,\dot{r})$ for all $u\in\mathbf{R}$, an improved Picard-Lindel\"{o}f theorem (see, e.g. \cite{Te})
	implies the existence of an even solution for all $u\in\mathbf{R}$. The solution is smooth since $f$ is smooth.
	Another application of the theorem to the original equation (\ref{DEP}) shows that this solution extends smoothly to the north and south poles of the sphere.

\vskip 8mm

\noindent
Zhongmin Shen \& Runzhong Zhao \\
Department of Mathematical Sciences \\
Indiana University-Purdue University Indianapolis \\
IN 46202-3216, USA  \\
E-mail: zshen@math.iupui.edu\\
E-mail: runzzhao@iu.edu

\end{document}